\documentclass{elsarticle}
\usepackage{amsmath,amssymb}
\usepackage{tikz}

\newtheorem{theorem}{Theorem}

\newtheorem{lemma}[theorem]{Lemma}
\newtheorem{corollary}[theorem]{Corollary}
\newtheorem{prop}[theorem]{Proposition}

\newproof{proof}{Proof}

\renewcommand\proof{\noindent\textsl{Proof. }}
\newcommand\sqr[2]{{\vbox{\hrule height.#2pt
    \hbox{\vrule width.#2pt height#1pt \kern#1pt
        \vrule width.#2pt}\hrule height.#2pt}}}
\renewcommand\qed{%
	\ifmmode\eqno\sqr53
	\else\nolinebreak\ \hfill\sqr53\medbreak\fi}
\tikzset
{
    treenode/.style = {circle, draw=black, align=center, minimum size=1cm},
    smalltree/.style = {circle, draw=black, align=center, minimum size=0.5cm},
    subtree/.style  = {circle,  align=center, minimum size =.5cm},
    phantom/.style = {circle}
}

\begin{document}

\title{Graphs that have a weighted adjacency matrix with spectrum
  $\{\lambda_1^{n-2}, \lambda_1^2\}$}

\author[karen]{Karen Meagher\corref{cor1}\fnref{fn1}}
\ead{karen.meagher@uregina.ca}

\address[karen]{Karen Meagher, Department of Mathematics and Statistics,\\
University of Regina, 3737 Wascana Parkway, Regina, SK, S4S 0A4, Canada}

\author[Irene]{Irene Sciriha}
\ead{irene.sciriha-aquilina@um.edu.mt}
\address[Irene]{Irene Sciriha, Department of Mathematics, \\ 
                        University of Malta, Msida, MSD2080, Malta}

\cortext[cor1]{Corresponding author}
\fntext[fn1]{Research supported by NSERC.}

\begin{abstract}
  In this paper we completely characterize the graphs which have an
  edge weighted adjacency matrix belonging to the class of $n \times
  n$ involutions with spectrum equal to $\{ \lambda_1^{n-2},
  \lambda_2^{2} \}$ for some $\lambda_1$ and some $\lambda_2$.  The
  connected graphs turn out to be the cographs constructed as the join of at least
  two unions of pairs of complete graphs, and possibly joined with
  one other complete graph.
\end{abstract}


\maketitle 
\section{Introduction}

To a graph $X$, we associate the collection of real $n \times n$
symmetric matrices defined by
\[ 
S(X) = \{ A \; : \; A=A^{T} ; \textrm{ for }  i \neq j, \;  a_{ij} \neq
0 \textrm{ if and only if } \{i,j\} \in E(X)\}.
\]
Note that there are no restrictions on the entries on the main
diagonal of a matrix in $S(X)$.  If $A \in S(X)$ for some graph $X$,
then $X$ is the \textsl{graph of $A$}.  This family of matrices has
been studied by many researchers and it is interesting to connect
properties of $S(X)$ to properties of the graph.  For example, there
has been significant work on determining the value of the minimum rank
over all matrices in $S(X)$ for a given graph $X$, see~\cite{psd, FH,
  FH2} and the references within.

For a real symmetric matrix $A$, let $q(A)$ denote the number of
distinct eigenvalues of $A$. For a graph $X$, define
\[
q(X) = \min\{q(A) \, : \, A \in S(X)\}.
\]
We say that $q(X)$ is the number of distinct eigenvalues for the graph
$X$. There have been several recent results regarding this parameter
\cite{MR3118943, BF, F, KS2, RUofWy}; this paper continues the work
in~\cite{MR3118943}. 

It is easy to see that $q(X) =1$ if and only if $X$ is an empty
graph. At the other extreme, $q(X) =|V(X)|$ if and only if $X$ is a
path~\cite{Fiedler}. There are very few known lower bounds on the
value of $q(X)$ for a graph $X$. One of the most effective is the
following simple bound; this is Theorem 3.2 from \cite{MR3118943}.

\begin{theorem}\label{thm:uniquepath}
Let $x$ and $y$ be two vertices of a graph $X$ at distance $d$.
If the path of length $d$ from $x$ to $y$ is unique, then $q(X) \geq  d + 1$.
\end{theorem}

The complete graph on $n$ vertices, denoted by $K_n$, has $q(K_n) =2$
(this can be achieved by the $(0,1)$-adjacency matrix of $K_n$).  Also
the complete bipartite graph $K_{n,n}$, and the hypercube also have
only two distinct eigenvalues.  The next example will show that the
family of graphs with only two distinct eigenvalues is very large.
For graphs $X$ and $Y$, the \textsl{join} of $X$ and $Y$, denoted $X
\triangledown Y$, is the graph with vertices $V(X) \cup V(Y)$, and
edge set
\[
E(X) \cup E(Y) \cup \{ \{x,y\} :  x \in V(X), \, y \in V(Y)\}. 
\]
In~\cite{MR3118943} it is shown that if $X$ is any connected graph,
then $q(X \triangledown X) =2$. From these examples,
it seems unlikely that the family of graphs $X$ with $q(X)=2$ can be
characterized.

For any graph, the multiplicities of the eigenvalues form an integer
partition of the number of vertices in the graph. In this paper, we
only consider graphs in which this partition has only two parts, so is
a \textsl{bipartition}. We say that $[n-i,i]$ is a
\textsl{multiplicity bipartition} of $X$, if there exists an $A \in
S(X)$ with spectrum $\{ \lambda_1^{n-i}, \lambda_2^i \}$.  The
\textsl{minimal} multiplicity bipartition of a graph is $[n-i,i]$ if
$i$ is the least value such that a multiplicity bipartition $[n-i,i]$
of the graph exists.  Note that if a non-empty graph $X$ has a
multiplicity bipartition, then $q(X) =2$. The main theorem of this
paper is a characterization of the graphs that have $[n-2,2]$ as their
minimal multiplicity bipartition.

\begin{theorem}\label{thm:main}
  Assume that $X$ is a connected graph. The minimal multiplicity
  bipartition of $X$ is $[n-2,2]$ if and only if
\[
X = (K_{a_1} \cup K_{b_1}) \triangledown (K_{a_2} \cup K_{b_2}) \triangledown \cdots \triangledown (K_{a_k} \cup K_{b_k})
\]
where $\{a_1,\dots,a_n\}$ and $\{b_1, \dots,b_n\}$ are non-negative
integers, $k>1$, and $X$ is not isomorphic to a complete graph, or to
$(K_{a_1} \cup K_{b_1}) \triangledown K_1$.
\end{theorem}

In the next section we shall state some known results related to graphs
with only two distinct eigenvalues. Section~\ref{sec:cographs}
presents some restrictions on graphs that have $[n-2,2]$ as a
multiplicity bipartition. In Section~\ref{sec:constr} we give
constructions for matrices $A \in S(X)$ that have spectrum
$\{\lambda_1^{n-2}, \lambda_2^2\}$ for all graphs $X$ identified in
Theorem~\ref{thm:main}. Section~\ref{sec:proofmain} gives the proof of
Theorem~\ref{thm:main}.

\section{Graphs with a multiplicity bipartition}
\label{sec:bipartition}

In this section, we present a way to determine if a
matrix has only two distinct eigenvalues with specific multiplicities.
The first result follows from \cite[Lemma 2.3]{MR3118943}.

\begin{prop}\label{prop:square2}
  Let $X$ be a non-empty graph, then $q(X) =2$ if and only if there is
  an $A \in S(X)$ with $A^2 = I$.
\end{prop}

We do not give a proof of Proposition~\ref{prop:square2}. Rather we
shall prove a stronger result, namely Lemma~\ref{lem:q2iffon}, which
implies it.  But first we need to introduce some notation.  Throughout
this paper, the $i$-th entry of a vector $v$ will be denoted by
$v_i$. We shall use $v(i)$ to denote different vectors (in
Section~\ref{sec:constr} we shall give constructions for vectors that
are based on a parameter $i$). The $j$-th entry of the vector $v(i)$
will be denoted by $v(i)_j$.  We shall start with a simple theorem
about the spectrum of a matrix with a specific form.

\begin{lemma}\label{lem:vectorconst}
Let $\{v(1), v(2), \dots, v(k)\}$ be a set of orthonormal vectors in
$\mathbb{R}^n$ with $1 \leq k <n$, and define
\[
A = I - 2(v(1)v(1)^T+v(2)v(2)^T+\dots +v(k)v(k)^T ).
\]
Then the following hold:
\begin{enumerate}
\item $A$ is real and symmetric;
\item the $(i,j)$-entry of $A$, where $i \neq j$, is
\[
-2 \sum_{\ell=1}^n u(\ell)_i u(\ell)_j;
\] 
\item the $(i,i)$-entry of $A$ is
\[
1 -2 \sum_{\ell=1}^n u(\ell)_i u(\ell)_i;
\] 
\item $A^2 = I$;
\item the spectrum of $A$ is $\{ 1^{(n-k)}, -1^{(k)}\}$; so $q(A) =2$; \label{cond5}
\item if $X$ is the graph of $A$, then $q(X) =2$. \label{cond6}
\end{enumerate}
\end{lemma}
\proof The first four statements follow immediately from the
definition of $A$.  The fifth follows from the fact that the set
$\{v(1),v(2), \dots,v(k)\}$ forms a set of orthogonal eigenvectors of
$A$ each with eigenvalue $-1$. Any vector orthogonal to all of $v(1),
v(2),\dots,v(k)$ is also an eigenvector of $A$, but with eigenvalue
$1$.  The final statement follows from Statement~\ref{cond5}, and the
fact that $X$ is non-empty.\qed

It is necessary that $k<n$ in the previous lemma, since if $k=n$, then
$A = -I$ and $q(A)=1$. The results that follow next show how the
previous lemma can be used to determine if a graph has a multiplicity
bipartition.

\begin{lemma}\label{lem:q2iffon}
Let $X$ be a non-empty graph. There exists $A \in S(X)$ with 
\[
A = I - 2(v(1)v(1)^T+v(2)v(2)^T+\dots +v(k)v(k)^T),
\]
where $\{v(1), v(2),\dots , v(k)\}$ is an orthonormal set of vectors
(with $1 \leq k<n$), if and only if $q(X) = 2$.
\end{lemma}
\proof 
By Statement~\ref{cond6} of Lemma~\ref{lem:vectorconst}, the condition
$q(x)=2$ is necessary. 

To prove that this condition is sufficient, let $B \in S(X)$ with $q(B) =2$.
Let $\lambda_1$ and $\lambda_2$ be the eigenvalues of $B$. Set
\begin{eqnarray}\label{eq:shiftandscale}
A = \frac{2}{\lambda_1 - \lambda_2}B - \frac{\lambda_1 + \lambda_2}{\lambda_1 - \lambda_2}I.
\end{eqnarray}
Then $A \in S(X)$ and the eigenvalues of $A$ are $-1$ and $1$.

By spectral decomposition, $A = 1 P_1 + (-1) P_{-1}$, where $P_i$ is
the projection to the $i$-eigenspace.  Since $P_1 + P_{-1} = I$, we
have that
\[
A = (I - P_{-1}) -P_{-1}  = I - 2 P_{-1}. 
\]
If $\{v(1),v(2),\dots,v(k)\}$ is an orthonormal basis for the $-1$ eigenspace, then 
\[
P_{-1} = v(1)v(1)^T+v(2)v(2)^T+\dots +v(k)v(k)^T.
\]
Thus $A$ has the required form. \qed

Statement~\ref{cond5} of Lemma~\ref{lem:vectorconst} and 
the proof of Lemma~\ref{lem:q2iffon} provide a way to construct an
adjacency matris for any graph $X$ with $q(X) =2$.

\begin{corollary}\label{cor:niceform}
A graph $X$ has $[n-k,k]$ as a multiplicity bipartition if and only if there
exist an $A \in S(X)$ such that 
\[
A = I -2(v(1)v(1)^T+v(2)v(2)^T+\cdots +v(k)v(k)^T)
\]
where $\{v(1), v(2),\dots ,v(k) \}$ is an orthonormal set of vectors.
\end{corollary}

In Sections~\ref{sec:cographs} and~\ref{sec:constr}, we use this to
determine the graphs with $[n-2,2]$ as a multiplicity bipartition. But
first we use this give information about the structure of a graph $X$
with $q(X)=2$. A \textsl{coclique} in a graph is a set of vertices in
which no two are adjacent, a coclique is also known as an
\textsl{independent set}.

\begin{lemma}\label{lem:nococliques}
Let $\{v(1), v(2), \dots, v(k)\}$ be a set of orthonormal vectors in $\mathbb{R}^n$ and
\[
A = I - 2(v(1)v(1)^T+v(2)v(2)^T+\cdots +v(k)v(k)^T ).
\]
Let $X$ be the graph of $A$. Then, provided that $X$ does not contain
any isolated vertices, the graph of $A$ does not contain a coclique of size $k+1$.
\end{lemma}
\proof Let $X$ be the graph of $A$ and label the vertices in $X$ by
$1,\dots ,n$.  Assume that vertices $1,2,\dots,k+1$ form a coclique in
$X$. The $(i,j)$-entry of $A$ for $i \neq j$ is $-2 \sum_{\ell=1}^k
v(\ell)_i v(\ell)_j$. So for all $1 \leq i,j \leq k+1$
\begin{eqnarray}\label{eq:zeros}
\sum_{\ell=1}^k v(\ell)_i v(\ell)_j =0.
\end{eqnarray}

Consider the vectors $x(i) = (v(1)_i, v(2)_i, \dots , v(k)_i )$ for $i
= 1,2, \dots, k+1$.  The vector $x(i)$ cannot be the zero vector,
since that would imply the vertex $i$ is an isolated vertex in $X$.
From Equation~\ref{eq:zeros}, if $i \neq j$, then $x(i)$ is orthogonal
to $x(j)$. But this implies that there is a set of $k+1$ non-zero,
length-$k$ orthogonal vectors, which is clearly not possible. Hence
$X$ does not have a coclique with more than $k$ vertices. \qed
With this lemma, we can characterize the graphs that have
$[n-1,1]$ as a multiplicity bipartition.

\begin{corollary}\label{cor:q2rankn-1}
  The graph $X$ has $[n-1,1]$ as a multiplicity bipartition if and
  only if $X$ is a complete graph with isolated vertices.
\end{corollary}
\proof If $X$ has $[n-1,1]$ as a multiplicity bipartition, then, by
Corollary~\ref{cor:niceform}, there is a matrix $A \in S(X)$ with $A = I -
2uu^T$. From Lemma~\ref{lem:nococliques}, any two vertices in $X$ are
adjacent, unless one of them is isolated.

Conversely, let $X$ be the graph with a clique of size $n$ and $k$ isolated
points. Let $J_n$ be the $n\times n$ matrix and $I_n$ the $n \times n$
identity matrix and consider the matrix
\[
A = 
\left( \begin{array}{c|c}
J_n-I_n & 0 \\ \hline
0 & -I_k
\end{array} \right).
\]
The spectrum of $A$ is $\{n-1, -1^{n+k-1}\}$ and $A \in S(X)$. Since
$X$ is non-empty, it has multiplicity bipartition $[n-1,1]$.\qed

In the previous corollary, we found that a graph with two distinct
eigenvalues and one with multiplicity $1$ must be the complete graph
with isolated vertices. This completely characterizes the graphs with
minimal multiplicity bipartition $[n-1,1]$. In the sequel, we shall give a characterization of the
graphs with minimal multiplicity bipartition $[n-2,2]$. In other words, we
characterize the graphs $X$ for which there is an $A \in S(X)$ with
\[
A = I -2(uu^T + vv^T)
\]
where $u$ and $v$ are orthonormal vectors.  If $u = (u_1,u_2, \dots,
u_n)$ and $v=(v_1,v_2, \dots, v_n)$, then the $(i,j)$-entry of $A$ is
$-2(u_iu_j + v_iv_j)$. So we need to determine which zero/non-zero
patterns are possible in $A$, with the conditions that $u$ and $v$ are
orthogonal and normalized.

The first restriction on the graphs with $[n-2,2]$ as a multiplicity
bipartition follows directly from Lemma~\ref{lem:nococliques}.

\begin{corollary}\label{cor:no3coclique}
  A connected graph $X$ with $[n-2,2]$ as a multiplicity bipartition
  does not have a coclique of size $3$.
\end{corollary}

\section{Cographs}
\label{sec:cographs}

In this section, we give a major restriction on the structure of the
graphs that have $[n-2,2]$ as a multiplicity bipartition.

\begin{lemma}\label{lem:noP4}
A graph $X$ with $[n-2,2]$ as a multiplicity bipartition does not have
an induced path of length three.
\end{lemma}
\proof Assume that the first four vertices in $X$ (which we simply
label $1$, $2$, $3$, $4$) form an induced path of length-$3$ in $X$. Assume
that there is a $A \in S(X)$ with
\[
A = I -2 (uu^T + vv^T).
\]
Then we have the following six equations:
\begin{equation}\label{eq1}
 \begin{aligned}
u_1u_3 + v_1v_3 = 0, \quad  &u_1u_4 + v_1v_4 = 0, \quad &u_2u_4+ v_2v_4= 0, \\
u_1u_2 + v_1v_2 \neq 0, \quad  &u_2u_3 + v_2v_3 \neq 0, \quad &u_3u_4 + v_3v_4 \neq 0. 
 \end{aligned}
\end{equation}

First we shall show that none of $u_i$ for $i \in \{1,\dots,4\}$ can
be equal to zero.  First suppose $u_3 = 0$. Then $v_1v_3 = 0$. If $v_3
= 0$, then the third vertex is not adjacent to the second vertex
(equation $u_2u_3 + v_2v_3 \neq 0 $ cannot hold). So $v_1 = 0$. But
then, since $u_1u_4+v_1v_4=0$, one of $u_1$ or $u_4$ is zero. If $u_1
= 0$, then the first vertex is not adjacent to the second vertex
(equation $u_1u_2 + v_1v_2 \neq 0 $ cannot hold). So this implies that
$u_4 =0$.  Then the equation $u_2u_4+ v_2v_4= 0$ implies that
$v_2v_4=0$.  If $v_2 =0$, then $u_2u_3 + v_2v_3 = 0$, which is a
contradiction. Similarly, if $v_4=0$, then $u_3u_4 + v_3v_4 = 0$,
which is also a contradiction. Thus $u_3 \neq 0 $. 

Similarly, we can show that $u_1,u_2$ and $u_4$ are also non-zero and
also that the entries in $v_i$ are not zero for $i \in \{1,\dots,4\}$.

Now we can assume that $u_i$ and $v_i$ are not
zero for $i \in \{1,\dots,4\}$. With this assumption, From
Equation~(\ref{eq1}), we have that
\[
u_1 = -\frac{v_1v_3}{u_3} = -\frac{v_1v_4}{u_4}.
\]
We set $k= \frac{v_3}{u_3} =
\frac{v_4}{u_4}$.  Similarly, 
\[
u_4 = -\frac{v_1v_4}{u_1} = -\frac{v_2v_4}{u_2},
\]
and in this case we set $\ell = \frac{v_1}{u_1} =
\frac{v_2}{u_2}$. Thus we have that
\[
v_1 = \ell u_1, \quad  v_2 = \ell u_2, \quad v_3 = k u_3, \quad v_4 = k u_4.
\]

Since vertices 1 and 4 are non-adjacent, 
\[
u_1u_4 + v_1v_4 = u_1u_4 + \ell u_1 k u_4 = 0,
\]
thus $1+\ell k = 0$.
But this implies that
\[
u_2u_3 + v_2v_3 = u_2u_3 + \ell u_2 k u_3 = (1+\ell k) u_2 u_3 = 0,
\]
which is a contradiction, since vertices $2$ and $3$ are adjacent.\qed

The family of graphs that do not contain a copy of $P_4$ are known as
the \textsl{cographs}. Cographs are a well-studied family of
graphs~\cite{MR1686154, MR1417258}. These graphs can be built
recursively.

\begin{prop}
The following recursive construction defines all cographs:
\begin{enumerate}
\item a single vertex is a cograph;
\item the union or join of two cographs is again a cograph; and
\item if $X$ is a cograph, then $\overline{X}$ is itself a cograph.
\end{enumerate}
\end{prop}

A \textsl{cotree} is a tree that is used to represent a cograph. There
is a 1-1 correspondence between cotrees and cographs; a cograph has a
unique cotree, and each cotree determines a unique cograph.The leaves
in the cotree correspond to vertices in the cograph. Internal nodes of
a cotree are labeled with either a union or a join. The children of
the nodes are connected by the operation by which the node is
labelled.

Let $T$ be the cotree of a cograph $X$.  We can assume that
below any internal node there must be at least two children, since if
there is just one, then the branch can be shortened.  Each child
represents a subgraph of $T$, that corresponds to a subgraph of $X$
that is also a cograph.  If an internal node is labeled as union, then
its children are either leaves, or internal nodes labeled with a
join. Similarly, if an internal node is labeled with a join, then its
children are either leaves or internal nodes labeled as unions.

We shall first characterize the cographs that do not contain a
coclique of size three. This will give a considerable restriction on
the possible graphs that have $[n-2,2]$ as a multiplicity bipartition.

\begin{lemma}\label{lem:cono3coclique}
If $X$ is a connected cograph with no coclique of size $3$, then 
\[
X = (K_{a_1} \cup K_{b_1}) \triangledown (K_{a_2} \cup K_{b_2}) \triangledown \cdots \triangledown (K_{a_k} \cup K_{b_k})
\]
where $\{a_1,\dots,a_k\}$ and $\{b_1 \dots,b_k\}$ are non-negative
integers.
\end{lemma}
\proof Assume that $X$ is a cograph with no cocliques of size three
and $T$ is the cotree for $X$.  Since $X$ is connected, the root of
$T$ must be a join.

If a vertex of $T$ that is labeled with a union has three children, then any set
formed by taking one vertex from the subgraph corresponding to each of
the children will be a coclique in $X$ of size three. Thus, in a
cograph with no coclique of size three in $X$, any internal vertex
labeled with a union can have at most two children, and each child
must correspond to a clique in $X$.

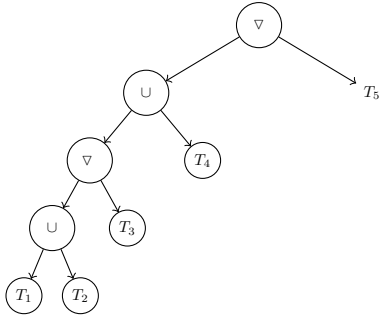
\begin{figure}[t]
\centering
{
    \begin{tikzpicture}[->,level/.style={sibling distance = 5cm/#1,
    level distance = 1.5cm},scale=0.6, transform shape]
    \node [treenode] {$\triangledown$}
    child
    {
        node [treenode] {$\cup$} 
        child
        {
            node [treenode] {$\triangledown$ } 
            child
            {
                node [treenode] {$\cup$} 
                        child{node [smalltree] {$T_1$} }
                        child{ node[smalltree] {$T_2$} }  
            }
            child
            {
                node [smalltree] {$T_3$} 
            }
        }
        child
        {
            node [smalltree] {$T_4$} 
        }
    }
    child
    {
        node [subtree] {$T_5$} 
    }
;
    \end{tikzpicture}
}
\caption{A cotree with an internal vertex labeled with a union with a
  descendent also labeled with a union.\label{fig:cotree}}
\end{figure}

If a vertex in $T$ is labeled with a union, then it cannot have a
descendent that is also labeled with a union. To see this, consider
Figure~\ref{fig:cotree}. This is an example of a cotree in which a
vertex labeled with a union has a descendant that is also labeled by a
union. Any set of three vertices in which one vertex is from each of the subgraphs
of $X$ corresponding to the cotrees $T_1, T_2, T_4$, will form a coclique of
size three in the graph $X$.

From these facts, the result holds.\qed

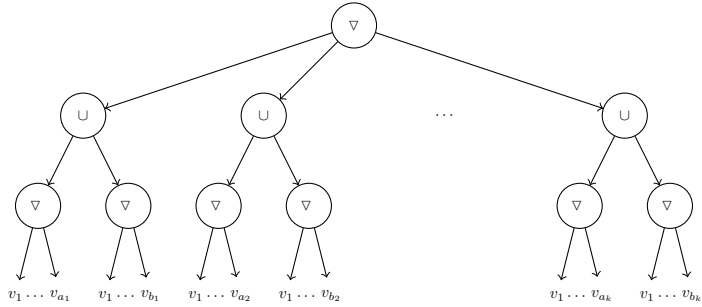
\begin{figure}[h!]
\centering
{
    \begin{tikzpicture}[->,level/.style={sibling distance = 10cm/(1.5*#1),
    level distance = 2cm},scale=0.6, transform shape]
    \tikzstyle{level 1}=[sibling distance=40mm] 
\tikzstyle{level 2}=[sibling distance=20mm] 
 \tikzstyle{level 3}=[sibling distance=10mm] 

\node [treenode] {$\triangledown$}
  child {node [treenode] (union1) {$\cup$} 
              child {
                 node [treenode](join1) {$\triangledown$ } 
%                    child{ node [subtree](v10) {$v_1$} }
                    child{ node [subtree](v11) {$v_1$} }
                    child{ node [subtree] (v12) {$v_{a_1}$} }
                   }
              child {
                 node [treenode](join1) {$\triangledown$ } 
%                    child{ node [subtree](v20) {$v_1$} }
                   child{ node [subtree](v21) {$v_1$} }
                    child{ node [subtree] (v22) {$v_{b_1}$} }
                    }
    }
 child {node [treenode] (union2) {$\cup$} 
              child {
                 node [treenode](join1) {$\triangledown$ } 
%                    child{ node [subtree](v30) {$v_1$} }
                   child{ node [subtree](v31) {$v_1$} }
                    child{ node [subtree] (v32) {$v_{a_2}$} }
                    }
              child {
                 node [treenode](join1) {$\triangledown$ } 
%                    child{ node [subtree](v40) {$v_1$} }
                   child{ node [subtree](v41) {$v_1$} }
                    child{ node [subtree] (v42) {$v_{b_2}$} }
                    }
    }
child [dashed, white] {node [phantom] (name) {}}
  child {node [treenode] (union3) {$\cup$} 
              child {
                 node [treenode](join1) {$\triangledown$ } 
 %                   child{ node [subtree](v50) {$v_1$} }
                    child{ node [subtree](v51) {$v_1$} }
                    child{ node [subtree] (v52) {$v_{a_k}$} }
                    }
              child {
                 node [treenode](join1) {$\triangledown$ } 
%                    child{ node [subtree](v60) {$v_1$} }
                   child{ node [subtree](v61) {$v_1$} }
                    child{ node [subtree] (v62) {$v_{b_k}$} }
                    }
    }
;
\path (union2) -- (union3) node [midway] {${\Large \cdots}$} ;
\path (v11) -- (v12) node [midway] {$\cdots$} ;
\path (v21) -- (v22) node [midway] {$\cdots$} ;
\path (v31) -- (v32) node [midway] {$\cdots$} ;
\path (v41) -- (v42) node [midway] {$\cdots$} ;
\path (v51) -- (v52) node [midway] {$\cdots$} ;
\path (v61) -- (v62) node [midway] {$\cdots$} ;
    \end{tikzpicture}
}
\caption{Cotrees for the cographs with no coclique of size three\label{fig:goodcotree}}
\end{figure}

In the following sections we shall consider cographs of the form
\[
X = (K_{a_1} \cup K_{b_1}) \triangledown (K_{a_2} \cup K_{b_2}) \triangledown \cdots \triangledown (K_{a_k} \cup K_{b_k}).
\]
We shall refer to a subgraph $K_{a_i} \cup K_{b_i}$ as a
\textsl{block} of $X$. Figure~\ref{fig:goodcotree} shows the general
form of a cotree for any such cograph.

\section{Constructions}
\label{sec:constr}

In this section, we shall give several results that are of the same
style.  In each result, for a given a graph $X$ of a certain form, we
construct a matrix $A \in S(X)$ with $A =I - 2( uu^T + vv^T)$ and $A^2
= I$. Note the condition from Lemma~\ref{lem:q2iffon} that $u$ and $v$
both have norm $1$, can be replaced with the condition that both
vectors have the same norm; in this case we use the matrix $A = I- (
2/\|u\|^2) (uu^T+vv^T)$.  Each result in this section gives examples
of two vectors $u$ and $v$ that they satisfy the following three
conditions:
\begin{enumerate}
\item $u$ and $v$ are orthogonal;
\item $u$ and $v$ have the same norm;
\item $u_iu_j +v_iv_j$ is zero if and only if vertices $i$ and $j$ in
  $X$ are non-adjacent.
\end{enumerate}

\begin{prop}\label{prop:nonemptyunions}
Let $a_i$ and $b_i$ be positive integers.
Suppose 
\[
X = (K_{a_1} \cup K_{b_1}) \triangledown (K_{a_2} \cup K_{b_2})
\triangledown \cdots \triangledown (K_{a_k} \cup K_{b_k}),
\]
then $X$ has $[n-2,2]$ as a multiplicity bipartition.
\end{prop}
\proof  
For $i = 1, \dots, k$, assume that $0 < a_o \leq b_i$.

For each subgraph $K_{a_i} \cup K_{b_i}$ of $X$, consider the two vectors
\begin{align*}
u(i)&=  \big( \underbrace{1, 1, \dots,  1,}_{a_i \textrm{ times }} \;
  \underbrace{ -\sqrt{\frac{a_i}{b_i}} w_i ,-\sqrt{\frac{a_i}{b_i}} w_i , \dots, -\sqrt{\frac{a_i}{b_i}} w_i}_{b_i \textrm{ times }}  \big) \\
v(i)&=\big( \underbrace{ w_i, w_i, \dots w_i,}_{a_i \textrm{ times}} \; 
        \underbrace{ \sqrt{\frac{a_i}{b_i}},\sqrt{\frac{a_i}{b_i}},
          \dots, \sqrt{\frac{a_i}{b_i}}}_{b_i \textrm{ times }}\big)
\end{align*}
where $w_i$ is any non-zero number. (The vertices of $K_{a_i} \cup
K_{b_i}$ are sorted so that the vertices in $K_{a_i}$ are first.)  The
norm of $u(i)$ equals the norm of $v(i)$, and the two vectors are
orthogonal.

Let $u$ be the vector formed by concatenating the vectors $u(1), u(2),
\dots, u(k)$, and $v$ be the vector formed from concatenating the vectors
$v(1), v(2), \dots, v(k)$. Then $u$ and $v$ have the same norm and are
orthogonal.

Finally we need to show that $A = I-(2/\|u\|^2)(uu^T + vv^T)$ is in
$S(X)$. 

If two vertices $x,y$ are from the same block in this graph,
say $K_{a_i} \cup K_{b_i}$, then
\[
[A]_{x,y} =
\begin{cases}
1+w_i^2, & \textrm{ if } x, y \in V(K_{a_i}); \\
0, & \textrm{ if } x \in V(K_{a_i}) \textrm{ and }  y \in V(K_{b_i}); \\
\frac{a_i}{b_i} (1+w_i^2), & \textrm{ if } x, y \in V(K_{b_i}).
\end{cases}
\]
So for any $x,y$ that are vertices in the block $K_{a_i} \cup
K_{b_i}$, the $(x,y)$-entry of $A$ will be zero if and only if $x$ and
$y$ are not adjacent in $X$.

Next consider the two vertices $x,y$ are from different blocks in this graph.
Assume $x$ is from $K_{a_i} \cup K_{b_i}$ and $y$ from $K_{a_j} \cup
K_{b_j}$, where $i \neq j$, then the corresponding entry in $A$ is
\[
[A]_{x,y} =
\begin{cases}
1+w_iw_j, & \textrm{ if } x \in V(K_{a_i}),\textrm{ and } y \in V(K_{a_j}); \\
-\sqrt{\frac{a_j}{b_j}} w_j + \sqrt{\frac{a_j}{b_j}}w_i ,
                & \textrm{ if }  x \in V(K_{a_i}) \textrm{ and } y \in V(K_{b_j}); \\
-\sqrt{\frac{a_i}{b_i}}w_i + \sqrt{\frac{a_i}{b_i}} w_j,
                & \textrm{ if } x \in V(K_{b_i})\textrm{ and } y \in V(K_{a_j}); \\
\sqrt{\frac{a_ia_j}{b_ib_j}} w_iw_j +  \sqrt{\frac{a_ia_j}{b_ib_j}}, 
               & \textrm{ if } x \in V(K_{b_i}) \textrm{ and } y \in V(K_{b_j}).
\end{cases}
\]
These reduce to 
\begin{align}\label{eq:fourvalues}
[A]_{x,y} =
\begin{cases}
1+w_iw_j,                                 & \textrm{ if } x \in V(K_{a_i}),\textrm{ and } y \in V(K_{a_j}); \\
\sqrt{\frac{a_j}{b_j}}(w_i -w_j),  & \textrm{ if } x \in V(K_{a_i}) \textrm{ and } y \in V(K_{b_j}); \\
\sqrt{\frac{a_i}{b_i}}(w_j-w_i) ,  &\textrm{ if }  x \in V(K_{b_i})\textrm{ and } y \in V(K_{a_j}); \\
\sqrt{\frac{a_ia_j}{b_ib_j}} (w_iw_j + 1), &\textrm{ if }  x \in V(K_{b_i}) \textrm{ and } y \in V(K_{b_j}).
\end{cases}
\end{align}
Since any two vertices from different blocks in $X$ are adjacent, we
need that the four values in Equation~\ref{eq:fourvalues} are all
non-zero.  It suffices to choose the $w_i$ so that $w_iw_j \neq -1$
and $w_i \neq w_j$. If we simply let $w_1, w_2, \dots, w_k$ be any set
of distinct positive real numbers, then the vectors $u$ and $v$ satisfy the
conditions.  The result follows from Corollary~\ref{cor:niceform}.\qed

In the previous result we required that all of the $a_i$ and $b_i$ be
larger than $0$. Next we shall consider the cases in which this
condition is dropped. If $a_i=b_i=0$ then $K_{a_i} \cup K_{b_i}$ has
no vertices, so we do not include it in $X$. If $a_i=0$, then $K_{a_i}
\cup K_{b_i} = K_{b_i}$ and including it is equivalent to the join
with a complete graph. In fact, if several of the $a_i = 0$, the
result is equivalent to the join with one large complete graph.  We
shall show that we can join a complete graph on at least two vertices
to the graphs in Proposition~\ref{prop:nonemptyunions} without
changing the multiplicity bipartition of the graph.

\begin{prop}\label{prop:singleemptyunions}
  Let $a_i$ and $b_i$ be positive integers and $c\geq 2$.  Suppose 
\[
X = (K_{a_1} \cup K_{b_1}) \triangledown (K_{a_2} \cup K_{b_2})
\triangledown \cdots \triangledown (K_{a_k} \cup K_{b_k}) \triangledown K_c,
\] 
then $X$ has $[n-2,2]$ as a multiplicity bipartition.
\end{prop}
\proof We shall construct the vectors $u$ and $v$ by concatenating
vectors $u(i)$ and $v(i)$ defined on the subgraphs $K_{a_i} \cup
K_{b_i}$, and vectors $u(k+1)$ and $v(k+1)$ defined on the subgraph
$K_{c}$.  For each subgraph $K_{a_i} \cup K_{b_i}$ in $X$, use the
vectors $u(i)$ and $v(i)$ defined in
Proposition~\ref{prop:nonemptyunions}. Next we define the vectors
$u(k+1)$ and $v(k+1)$ which are indexed by the vertices in the
subgraph $K_c$.

If $c=2$ simply set these two vectors to be
\[
u(k+1)=(1,0), \quad v(k+1)=(0,1).
\]
Since the entries of $u$ and $v$ are all non-zero, there will be edges
between any vertex in $K_2$ and any vertex not in $K_2$. 

If $c \geq 3$, define the following two vectors
\begin{align*}
u(k+1)&=\big( -1, 1, \dots , 1, -\frac{c-2}{2} \big),\\
v(k+1)&=\big( \frac{c-2}{2}, 1, 1, \dots, 1 \big).
\end{align*}

In Proposition~\ref{prop:nonemptyunions}, the values of the $w_i$ could
be any distinct positive real numbers. For the vectors $u$ and $v$
defined here to  satisfy the condition that any vertex from $K_c$ is
adjacent to any vertex in $K_{a_i} \cup K_{b_i}$, it is necessary that
the following six conditions hold:
\begin{enumerate}
\item $(-1)(1)+ \left(\frac{c-2}{2}\right)(w_i) \neq 0$ ,
\item $(1)(1) + (1)(w_i) \neq 0$ ,
\item $(-\frac{c-2}{2})(1) + (1)(w_i) \neq 0$,
\item $(-1)\left(-\sqrt{\frac{a_i}{b_i}} w_i \right) + \left(\frac{c-2}{2}\right)\left(\sqrt{\frac{a_i}{b_i}} \right)\neq 0$,
\item $(1) \left(-\sqrt{\frac{a_i}{b_i}} w_i \right) + (1) \left
    (\sqrt{\frac{a_i}{b_i}}\right) \neq 0$,
\item $\left(-\frac{c-2}{2}\right) \left(-\sqrt{\frac{a_i}{b_i}}w_i \right) +(1) \left(\sqrt{\frac{a_i}{b_i}} \right) \neq 0$.
\end{enumerate}
This requires simple additional restrictions on the values of
$w_i$. Specifically, the $w_i$ need to be positive, distinct integers
that are not in the set $\{\pm 1, \pm \frac{c-2}{2}, \pm
\frac{2}{c-2}\}$. Again, by Corollary~\ref{cor:niceform}, $X$ has
$[n-2,2]$ as a multiplicity bipartition.
\qed

The final case that we consider is 
\[
X = (K_{a_1} \cup K_{b_1}) \triangledown (K_{a_2} \cup K_{b_2}) \triangledown \cdots \triangledown
(K_{a_k} \cup K_{b_k}) \triangledown K_1.
\]
This construction will be based on the construction in
Proposition~\ref{prop:singleemptyunions}, but first we give a
construction for a subgraph.

\begin{prop}\label{prop:oneone}
Let $X = (K_a \cup K_b) \triangledown (K_c \cup K_d) \triangledown
K_1$, with $a,b,c,d$ all positive integers. Then $X$ has $[n-2,2]$ as a multiplicity bipartition.
\end{prop}
\proof
Assume that $0<a \leq b$ and $0 <c \leq d$.

For the subgraph $K_a \cup K_b$ define 
\begin{align}
u(1) = \big( \underbrace{1, \dots, 1,}_{a \textrm{ times}} \, 
                     \underbrace{4, \dots, 4}_{b \textrm{ times}}
                     \big), 
\qquad
v(1) = \big( \underbrace{2, \dots, 2,}_{a \textrm{ times}} \, 
                     \underbrace{-2, \dots,-2}_{b \textrm{ times} } \big)
\end{align}
and for the subgraph $K_c \cup K_d$ define
\begin{align}
u(2) = \big( \underbrace{2, \dots, 2,}_{c \textrm{ times}} \,
                     \underbrace{9, \dots, 9}_{d \textrm{ times}}
                     \big),
\qquad
v(2) = \big(\underbrace{6, \dots, 6}_{c \textrm{ times}} , \,
\underbrace{-3, \dots, -3}_{d \textrm{ times}} \big)
\end{align}
Define $u = (u(1), \, u(2), \, x )$ and $v=(v(1), \, v(2), \,
y)$. 

Then
\[
u\cdot v = 2a-8b+12c-27d +xy.
\]
Since $u$ and $v$ must be orthogonal, it is necessary that
\[
x = \frac{ -2a+8b-12c+27d}{y}.
\]
As $0< a \leq b$ and $0<c\leq d$, the numerator is strictly
positive; this implies that $x$ and $y$ have the same sign, and neither
are equal to $0$.

Further,
\begin{align*}
\|u\|^2 = 1a+16b+4c+81d + x^2, \qquad \|v\|^2 = 4a+4b+36c+9d+y^2.\\
\end{align*}
The vectors $u$ and $v$ must have the same norm, which implies that 
\[
x^2 = 3a -12b +32c-72d +y^2.
\] 
Eliminating $x$, this becomes
\[
\left( \frac{ -2a+8b-12c+27d}{y} \right)^2 = 3a -12b +32c-72d +y^2
\]
which gives 
\begin{align}
0 = y^4 + (3a -12b +32c-72d )  y^2 - (-2a+8b-12c+27d)^2.  \label{eq:quadratic}
\end{align}

Let $A= 3a -12b +32c-72d $ and $B = -2a+8b-12c+27d$. Consider the function 
\[
f(s) = s^2 + A s - B^2
\]
Then $f(2B) = 3B^2+2AB = B(28c-63d)$. Since $B$ is positive and $c\leq
d$, it follows that $f(2B)$ is strictly negative. Similarly, $f(3B) =
8B^2+3AB = B(-7a+28b)$, which, since $a\leq b$, is strictly positive.
Thus $f$ has at least one root $r$ in the range $(2B, 3B)$. 

Set $y = \sqrt{r}$ (since $r$ is positive, this is possible). Then $y$
satisfies Equation~\ref{eq:quadratic}. We also set $x
=(-2a+8b-12c+27d)/y$ (this implies that both $x$ and $y$ are
positive).  At this point, we have defined the entries of $u$ and $v$
so that they are orthogonal and have the same norm. Finally, we need
to show that $u_iu_j +v_iv_j = 0$ if and only if $i$ and $j$ represent
a pair of non-adjacent vertices.

If one of vertices $i$ and $j$ is from $K_a$ and the other from $K_b$, (or one
from $K_c$ and the other from $K_d$), then it is clear that $u_iu_j
+v_iv_j = 0$. Further, for any other pair of vertices, both from $(K_a
\cup K_b) \triangledown (K_c \cup K_d)$, the value of $u_iu_j +v_iv_j$
is not equal to zero.

Next consider the case where $i$ is a vertex in $K_1$ and $j$ is in $(K_a \cup K_b)
\triangledown (K_c \cup K_d)$, then $u_iu_j +v_iv_j $ is one of
\[
x+2y,\quad 4x-2y, \quad 2x+6y, \quad 9x-3y.
\]
Since $x$ and $y$ are both positive, $x+2y \neq 0$ and $2x+6y \neq 0$.
If $4x-2y=0$, then, since $xy = -2a+8b-12c+27d$, this
implies that $2 (-2a+8b-12c+27d) =y^2$. Similarly, if $9x-3y = 0$,
then $3(-2a+8b-12c+27d)=y^2$.  But, both of these are impossible,
since $y^2$ is equal to a root of $f(s)$ strictly between $2(
-2a+8b-12c+27d)$ and $3(-2a+8b-12c+27d)$. Thus, by
Corollary~\ref{cor:niceform}, $X$ has $[n-2,2]$ as a multiplicity bipartition.\qed

\begin{prop}\label{prop:completeandlots}
Let 
\[
X = (K_{a_1}\cup K_{b_1}) \triangledown (K_{a_2} \cup K_{b_2}) \triangledown
 \cdots \triangledown (K_{a_k} \cup K_{b_k}) \triangledown K_1,
\] 
where $k\geq 2$. Then $X$ has $[n-2,2]$ as a multiplicity bipartition.
\end{prop}
\proof Use the vectors $u(1)$, $v(1)$, $u(2)$ and $v(2)$ defined in
the proof of Proposition~\ref{prop:oneone} for the blocks $K_{a_1}\cup
K_{b_1}$ and $K_{a_2} \cup K_{b_2}$. Further, use the length-one
vectors $u(k+1) = (x)$ and $v(k+1) = (y)$ that are defined in
Proposition~\ref{prop:oneone} for the block $K_1$. For all other
$K_{a_i} \cup K_{b_i}$ use the vectors $u(i)$ and $v(i)$ defined in
the proof of Proposition~\ref{prop:nonemptyunions}. Define $u$ to be
the vector formed by concatenating $u(1), \dots, u(k+1)$ and $v$ the
vectors formed by concatenating $v(1), \dots, v(k+1)$.

Provided that the $w_i$ are distinct positive numbers not in the set
\[
\{ \pm 2, \pm 1/2,  \pm 1/3, \pm 3, \pm x/y, \pm y/x\},
\]
then $A = I-2(uu^T + vv^T) \in S(X)$ and by
Corollary~\ref{cor:niceform}, $X$ has $[n-2,2]$ as a multiplicity
bipartition.  \qed

\section{Proof of Main Theorem}
\label{sec:proofmain}

We now have all the tools to give the exact characterization of graphs
with minimal multiplicity bipartition $[n-2,2]$. 

{\bf Proof of Theorem~\ref{thm:main}.}
If $X$ has $[n-2,2]$ as a multiplicity bipartition, then by Lemma~\ref{lem:noP4}, it
is a cograph and by Corollary~\ref{cor:no3coclique} it has no
cocliques with three vertices. Then by Lemma~\ref{lem:cono3coclique},
the graph must have the form
\[
X = (K_{a_1} \cup K_{b_1}) \triangledown (K_{a_2} \cup K_{b_2}) \triangledown \cdots \triangledown (K_{a_k} \cup K_{b_k})
\]
(we assume that $b_i \geq a_i$ and that $(a_i ,b_i) \neq (0,0)$).

If $a_1 =0$ and $a_2=0$ then 
\[
X = (K_0 \cup K_{b_1 + b_2}) \triangledown \cdots \triangledown (K_{a_n} \cup K_{b_n})
\]
so we can assume that $a_i =0$ for at most one $i$.

Assume that there is one $i$ with $a_i =0$. If $b_i \geq 2$, then by
Proposition~\ref{prop:singleemptyunions} the graph has $[n-2,2]$ as a
multiplicity bipartition.  Provided that $X$ is not the complete graph
(using Corollary~\ref{cor:q2rankn-1}), this implies that the minimal
multiplicity bipartition of $X$ is $[n-2,2]$.

If $a_i=0$, $b_i =1$ and $k \geq 3$, then by
Proposition~\ref{prop:completeandlots} the graph has $[n-2,2]$ as a
multiplicity bipartition. The fact that $k\geq 3$ ensures that the
graph is not complete, so this is the minimal multiplicity
bipartition.

If $a_i=0$, $b_i=1$ and $k=2$, then $X$ is equal
to
\[
X =  (K_{a_1} \cup K_{b_1}) \triangledown K_1.
\]
In this case there is a unique path with two edges from any vertex in
$K_{a_1}$ to any vertex in $K_{b_1}$. By Theorem~\ref{thm:uniquepath}, this
implies that $q(X) \geq 3$ and that this graph has no multiplicity
bipartition. Furthermore, if $a_1=0$, $b_1=1$ and $k=1$, then $X$ is just
a single isolated vertex, in which case the spectrum contains just one
eigenvalue.

Finally, if all the $a_i$ and $b_i$ are greater than 0, then $X$ is not
a complete graph and has minimal multiplicity bipartition $[n-2,2]$ by
Proposition~\ref{prop:nonemptyunions}.  \qed

We can also characterize the disconnected graphs with minimal
multiplicity bipartition $[n-2,2]$.

\begin{theorem}\label{thm:submain}
A disconnected graph $X$ has minimal multiplicity bipartition $[n-2,2]$ if and only if
\[
X = K_{a} \cup K_{b} \cup \overline{K_c}, 
\]
or $X$ is the union of a graph with minimal multiplicity bipartition $[n-2,2]$ and isolated vertices.
\end{theorem}
\proof
Assume that $X = K_{a} \cup K_{b} \cup \overline{K_c}$. By
Corollary~\ref{cor:q2rankn-1}, there is a matrix $A \in S(K_a)$ with
spectrum $\{1^1,(-1)^{a-1}\}$, and $B \in S(K_{b} \cup
\overline{K_c})$ with spectrum $\{1^1,(-1)^{b+c-1}\}$. Then the matrix
\[
C = 
\left( \begin{array}{c|c}
A & 0 \\ \hline
0 & B
\end{array} \right).
\]
is in $S(X)$ and has spectrum $\{1^2,(-1)^{a+b+c-2}\}$.

Conversely, assume that there is an $A \in S(X)$ with spectrum
$\{\lambda_1^{n-2}, \lambda_2^2\}$. The spectrum for a component of
$X$ must be one of
\[
\{\lambda_1^{k-2}, \lambda_2^2\}, \quad  \{\lambda_1^{k-1}, \lambda_2^1\}, \quad \{\lambda_1^k\}.
\]

If one component has spectrum $\{\lambda_1^{k-2}, \lambda_2^2\}$, then
the spectrum of every other component must be $\{\lambda_1^{\ell_i}
\}$.  Thus the first component is a graph with multiplicity
bipartition $[n-2,2]$ and the remaining components are isolated
vertices.

If a component has spectrum $\{\lambda_1^{k-1}, \lambda_2^1\}$, then,
by Corollary~\ref{cor:q2rankn-1}, that component is a complete
graph. There can be at most two components with spectrum
$\{\lambda_1^{k-1}, \lambda_2^1\}$ and the remaining components must
be isolated vertices.  \qed

\section{Further Work}

In this paper we characterized the graphs that have minimal
multiplicity bipartition $[n-1,1]$ and $[n-2,2]$. An obvious next step
is to determine which graphs have $[n-k,k]$ as a multiplicity
bipartition for larger values of $k$. This may also be a route to
determining the entire family of graphs $X$ that have $q(X)=2$.

This problem can be generalized to any integer partition of $n$.  Let
$\pi = (\pi_1, \pi_2, \dots, \pi_k)$ be an integer
partition of $n$. We say that $\pi$ is a \textsl{multiplicity
  partition} of $X$, if there exists an $A \in S(X)$ with spectrum $\{
\xi^{\pi_1}, \xi^{\pi_2}, \dots , \xi^{\pi_k}\}$. In the
case where the partition has only two parts there is a natural concept
of the minimal partition. But, if $\pi$ is not a bipartition, it
is not clear what a minimal partition would be.

This concept relates both the minimal number of distinct eigenvalues
of a graph and the maximum multiplicity for a graph. For a graph $X$,
the value of $q(X)$ is equal to the fewest number of parts in a
multiplicity partition of the graph. The maximum multiplicity is the
largest size of a part in a multiplicity partition. 

An interesting example is the class of even cycles. The maximum
multiplicity of an eigenvalue for this graph is $2$~\cite{Fiedler} and
$q(C_{2k}) = k$~\cite[lemma2.7]{MR3118943}.  In~\cite[Thm. 3.3]{FF} it
is show that for any set of numbers $\lambda_1=\lambda_{2}>
\lambda_3=\lambda_{4}>\ldots > \lambda_{2k-1}=\lambda_{2k},$ there is
an $A\in S(C_{2k})$ with spectrum $\{\lambda_1, \lambda_2,\dots,\lambda_{2k}\}$.
This implies that $[2^{k}]$ is a multiplicity partition for
$C_{2k}$. Since no multiplicity partition can have fewer parts, nor
any parts of larger size, we claim that this is the minimal
multiplicity partition for the even cycles.

Another interesting family of graphs to consider are the paths. Let
$P_n$ denote the path on $n$ vertices. Every $A \in S(P_n)$ will have
$n$ distinct eigenvalues, so $[1^n]$ is a multiplicity partition of
$P_n$. In this case, the maximum multiplicity of the path is $1$ and
$q(P_n)=n$, so $[1^n]$ is the only multiplicity partition for
$P_n$.  In \cite[Section 7]{MR3118943} graphs $X$ with $q(X) =
|V(X)|-1$ are considered. These are the graph that have $[2,1^{n-2}]$,
(and possibly $[1^n]$) as a multiplicity partition, but no other
multiplicity partitions.

There are also many questions related to the existence of a spectrum
bipartition. For example, if a graph has $\lambda$ as a multiplicity
partition, does this imply that the graph also has spectrum partition
$\mu$ for some other partition $\mu$? Can these relations then be used
to define an ordering on the partitions?

\end{document}